\documentclass[12pt]{article}
\usepackage{a4wide}
\usepackage[centertags,leqno]{amsmath}
\usepackage{amssymb,amsfonts}
\usepackage{mathrsfs}
\usepackage{color}
\usepackage{graphicx}
\usepackage[colorlinks=true]{hyperref}

\newtheorem{theorem}{Theorem}[section]

\newtheorem{lemma}[theorem]{Lemma}
\newtheorem{proposition}[theorem]{Proposition}
\newtheorem{remark}[theorem]{Remark}

\definecolor{violet}{rgb}{0.5,0,0.5}
\definecolor{orange}{cmyk}{0,0.3,0.7,0}

\newtheorem{definition}[theorem]{Definition}

\numberwithin{equation}{section}

\newcommand{\qed}{\rule{2mm}{2mm}}

\newcommand{\eqdef}{\stackrel{{\mathrm {def}}}{=}}

\renewcommand{\colon}{:\,}
\newcommand{\proof}{{\em Proof. }}

\newcommand{\RR}{\mathbb{R}}

\begin{document}

\baselineskip=12pt
\begingroup\Large\bf
\begin{center}
  Travelling waves in the Fisher\--KPP equation\\
  with nonlinear diffusion and\\
  a non\--Lipschitzian reaction term
\end{center}
\endgroup

\vspace{0.3cm}
\begingroup\rm
\begin{center}
\centerline{Pavel {\sc Dr\'{a}bek}
}
\begin{tabular}{c}
  Department of Mathematics and\\
  N.T.I.S. (Center of New Technologies for Information Society)\\
  University of West Bohemia\\
  P.O.~Box 314,
  CZ-306 14 Plze\v{n}, Czech Republic\\
  {\it e}-mail: {\tt pdrabek@kma.zcu.cz}
\end{tabular}
\end{center}
\endgroup

\vspace{-0.5cm}
\begin{center}
and
\end{center}

\vspace{-0.0cm}
\begingroup\rm
\begin{center}
\centerline{Peter {\sc Tak\'{a}\v{c}}
}
\begin{tabular}{c}
  Institut f\"ur Mathematik,
  Universit\"at Rostock\\
  Ulmenstra{\ss}e~69, Haus~3,
  D-18055 Rostock, Germany\\
  {\it e}-mail: {\tt peter.takac@uni-rostock.de}
\end{tabular}
\end{center}
\endgroup

\vspace{-0.1cm}

\begin{center}
\today
\end{center}

\vspace{0.1cm}


\noindent
\begingroup\footnotesize
{\bf {\sc Abstract.}}
We consider a one\--dimensional reaction\--diffusion equation of
Fisher\--Kolmogoroff\--Petrovsky\--Piscounoff type.
We investigate the effect of the interaction between
the nonlinear diffusion coefficient and the reaction term
on the existence and non\-existence of travelling waves.
Our diffusion coefficient is allowed to be degenerate or singular
at both equilibrium points, $0$ and $1$, while
the reaction term need not be differentiable.
These facts influence the existence and qualitative properties of
travelling waves in a substantial way.
\endgroup

\vfil
\noindent
\begin{tabular}{ll}
{\bf Running head:}
& Travelling waves in the Fisher\--KPP equation\\
\end{tabular}

\noindent
\begin{tabular}{ll}
{\bf Keywords:}
& Fisher\--Kolmogoroff\--Petrovsky\--Piscounoff equation;
  travelling wave;\\
& degenerate and/or singular diffusion;
  non\-smooth reaction term;\\
& existence and non\-existence of travelling waves;\\
& an overdetermined first\--order boundary value problem.\\
\end{tabular}

\noindent
\begin{tabular}{ll}
{\bf 2010 Mathematics Subject Classification:}
& Primary   35Q92, 92D25, (34B08);\\
& Secondary 35K57, 35K65, (34B18).
\end{tabular}

\noindent
{\it Declarations of interest:\/} {\bf none}.
 
\newpage

\baselineskip=14pt
\section{Introduction}
\label{s:Intro}

We are concerned with the travelling waves
(particularly their speed and profile)
for the {\it Fisher\--Kolmogoroff\--Petrovsky\--Piscounoff\/}
population model with {\it nonlinear diffusion\/}
(of porous medium type)
and a {\it non\--Lipschitzian\/} reaction term,
\begin{equation}
\label{e:FKPP}
    \frac{\partial u}{\partial t}
  - \frac{\partial}{\partial x}
    \left( d(u)\, \frac{\partial u}{\partial x} \right)
  = g(u)
    \quad\mbox{ for }\, (x,t)\in \RR\times \RR_+ \,.
\end{equation}
We employ certain specific forms of
the possibly degenerate or singular diffusion coefficient $d(u)$
and the nonlinear reaction function $g(u)$
that are motivated by classical population models by
{\sc R.~A.\ Fisher} \cite{Fisher} and
{\sc A.~N.\ Kolmogoroff}, {\sc I.~G.\ Petrovsky}, and
{\sc N.~S.\ Piscounoff} \cite{KPP},
both from the same year of $1937$.
We allow both $d(u)$ and $g(u)$ to depend continuously
on the population density $u$.
The reaction\--diffusion equation \eqref{e:FKPP}
is briefly referred to as
the {\it Fisher\--KPP equation\/} (or {\it FKPP equation\/}).

In contrast with similar models that have been considered
in the literature so far
\cite{Malag-Marce-1} -- \cite{Sanchez-Maini},
our diffusion term $d = d(u)$ and the reaction term $g = g(u)$
are much more general functions.
In fact, the diffusion term $d = d(u)$ may degenerate or blow up as
$u\to 0+$ and/or $u\to 1-$.
At the same time, the reaction term $g = g(u)$ need not be
a Lipschitz continuous function in its domain of definition.
While the role of the nonlinear reaction term $g = g(u)$
has been justified already in the original works
\cite{Fisher, KPP}
(which consider only constant diffusivity $d > 0$),
the importance of the density\--dependent diffusion term $d = d(u)$
in insect despersal models is emphasized in
{\sc J.~D.\ Murray}'s monograph \cite[{\S}13.4, p.~449]{Murray-I}.

In a general biological {\it Fisher\--KPP model\/}
one naturally expects travelling waves $u(x,t) = U(x-ct)$
with a continuous wave profile $U$.
However, requiring a smoother profile $U$ does not seem to be
biologically justified, see
\cite[{\S}11.3]{Murray-I} for a sketch of non\-smooth profiles
in Fig.\ 11.2 on p.~403.
Taking into account this fact, we define a travelling wave for
problem \eqref{e:FKPP} in a rather general fashion
that does not require differentiability of the profile; cf.\
Definition~\ref{def-trav_wave} below.

Density\--dependent dispersal
(modelled by density\--dependent diffusion)
has been observed in many insect polulations, such as
the ant\--lion {\em Glenuroides japonicus\/}.
Several authors propose to analyse the flux of ants throughout
a compartmentally divided habitat which leads to
the spatial segregation of a species.
For greater details and numerous references to biological modelling,
we refer the reader to
{\sc F.\ S\'anchez\--Gardu\~{n}o} and {\sc P.~K.\ Maini}
\cite[Sect.~2, pp.\ 164--166]{Sanchez-Maini}.

This article is organized as follows.
Our new definition of a travelling wave is given in the next section
(Section~\ref{s:exist:FKPP}).
Basic properties of a wave profile $U$, such as monotonicity,
are studied in Section~\ref{s:Wave_Profile}.
A standard phase plane transformation applied to the equation
for the wave profile $U$ in Section~\ref{s:Phase_Plane}
yields an overdetermined first\--order, two\--point
boundary value problem.
This is our basic tool for obtaining existence and nonexistence of
a travelling wave.
The last section (Section~\ref{s:Trav-Shape})
is dedicated to studies with simple terms
$d(u)$ and $g(u)$ that are nonlinear of {\it power\--type\/}
near the equilibrium points.
As a conclusion, from the interaction between $d(u)$ and $g(u)$
we determine the asymptotic shape of travelling waves
near the equilibrium points.

\section{A quasilinear Fisher\--KPP equation\\
         with non\--smooth positive reaction}
\label{s:exist:FKPP}

As usual, we denote
$\RR\eqdef (-\infty, \infty)$, $\RR_+\eqdef [0,\infty)$,
and assume that the diffusion coefficient $d$ and the reaction term $g$
satisfy the following basic hypotheses, respectively:
\begin{enumerate}
\renewcommand{\labelenumi}{{\bf (H\alph{enumi})}}
\item[{\bf (H1)}]
$d\colon \RR\setminus \{0,\, 1\}\to \RR$
is a continuous, but
{\it not necessarily smooth\/} function, such that
$d(s) > 0$ for every $s\in \RR\setminus \{0,\, 1\}$, and
the (Lebesgue) integral
$\int_a^b d(s) \,\mathrm{d}s < \infty$ whenever
$-\infty < a < b < +\infty$.
\item[{\bf (H2)}]
$g\colon \RR\to \RR$ is a continuous, but
{\it not necessarily smooth\/} function,
such that
$g(0) = g(1) = 0$ together with
$g(s) > 0$ for every $s\in (0,1)$, and
$g(s) < 0$ for every $s\in (-\infty,0)\cup (1,\infty)$.
\end{enumerate}

\par\noindent
The reaction function $g$ satisfying {\bf (H2)} comprises also
the so\--called
{\em generalized logistic growth\/}
in the population model studied in
{\sc A.\ Tsoularis} and {\sc J.\ Wallace} \cite{Tsoular_Wall}.

We reformulate eq.~\eqref{e:FKPP} for $u(x,t)$
as an equivalent initial value problem for the unknown function
$v(z,t) = u(z+ct,t)\equiv u(x,t)$ with the moving coordinate $z = x-ct$:
\begin{equation}
\label{e_c:FKPP}
    \frac{\partial v}{\partial t}
  - \frac{\partial}{\partial z}
    \left( d(v)\, \frac{\partial v}{\partial z} \right)
  - c\, \frac{\partial v}{\partial z}
  = g(v) \,,
    \quad (z,t)\in \RR\times \RR_+ \,.
\end{equation}
We will show below that every travelling wave $u(x,t) = U(x-ct)$
for \eqref{e:FKPP} must have a monotone decreasing profile
$U\colon \RR\to \RR$ satisfying
\begin{equation}
\label{lim:FKPP}
  \lim_{z\to -\infty} U(z) = 1 \quad\mbox{ and }\quad
  \lim_{z\to +\infty} U(z) = 0 \,.
\end{equation}

More precisely, $U\colon \RR\to \RR$ must be monotone decreasing with
$U'< 0$ on a suitable open interval $(z_0,z_1)\subset \RR$, such that
\begin{equation}
\label{z_0,1:FKPP}
  \lim_{z\to z_0+} U(z) = 1 \quad\mbox{ and }\quad
  \lim_{z\to z_1-} U(z) = 0 \,,
\end{equation}
by Proposition~\ref{prop-monot_TW}.
We would like to remark that the cases of
$z_0 > -\infty$ and/or $z_1 < +\infty$
render qualitatively different travelling waves than
the classical case $(z_0,z_1) = \RR$ which has been studied
in the original works \cite{Fisher, KPP} and in the literature
\cite{AronWein-1, AronWein-2, Fife-McLeod, Hamel-Nadira,
      Malag-Marce-1, Murray, Murray-I, Sanchez-Maini}.

In order to be able to give a workable definition of a travelling wave,
we introduce the (Lebesgue) integral
\begin{equation*}
  D(s)\eqdef \int_0^s d(s') \,\mathrm{d}s'
  \quad\mbox{ for every }\, s\in \RR \,.
\end{equation*}
This is an absolutely continuous function on $\RR$
which is continuously differentiable on $\RR\setminus \{0,\, 1\}$
with the derivative $D'(s) = d(s)$ for every
$s\in \RR\setminus \{0,\, 1\}$.
Using this setting, in Section~\ref{s:Phase_Plane}
we are able to find a {\it first integral\/} for
the second\--order equation for $U$ restricted to the open interval
$(z_0,z_1)\subset \RR$:
\begin{equation}
\label{eq:FKPP}
    \frac{\mathrm{d}}{\mathrm{d}z}
    \left( d(U)\, \frac{\mathrm{d}U}{\mathrm{d}z} \right)
  + c\, \frac{\mathrm{d}U}{\mathrm{d}z} + g(U(z)) = 0\,,
    \quad z\in (z_0,z_1) \,.
\end{equation}
It is easy to observe that this equation is valid for every
$z\in \RR\setminus \{ z_0,\, z_1\}$
(in the sense of Definition~\ref{def-trav_wave} below)
provided $U$ is extended by
$U(z) = 1$ if $-\infty < z\leq z_0$ and
$U(z) = 0$ if $z_1\leq z < +\infty$.

\begin{definition}\label{def-trav_wave}\nopagebreak
\begingroup\rm
A function $u(x,t) = U(x-ct)$ of $(x,t)\in \RR\times \RR_+$
is called a {\it travelling wave\/}
(or {\em TW}, for short)
for problem \eqref{e:FKPP} where $c\in \RR$ is a constant called\/
{\it wave speed\/} (or simply {\it speed\/}) and\/
$U\colon \RR\to \RR$ is a continuous function called\/
{\it wave profile\/} (or simply {\it profile\/})
with the following properties:
\begin{itemize}
\item[{\rm (a)}]
$U(z)\geq 0$ holds for every $z\in \RR$ and
the limits in \eqref{lim:FKPP} are valid.
\item[{\rm (b)}]
The composition
$z\mapsto D(U(z))\colon \RR\to \RR$
is a continuously differentiable function with the derivative
\begin{math}
  \frac{\mathrm{d}}{\mathrm{d}z}\, D(U(z))
\end{math}
vanishing at every point\/ $\xi\in \RR$ such that\/
$U(\xi)\in \{ 0,\, 1\}$.
\item[{\rm (c)}]
The following integral form of eq.~\eqref{eq:FKPP}
is valid for all pairs\/ $z, z^{\ast}\in \RR$:
\begin{equation}
\label{eq:int_U(z)}
\begin{aligned}
&   \frac{\mathrm{d}}{\mathrm{d}z}\, D(U(z))
  - \frac{\mathrm{d}}{\mathrm{d}z}\, D(U(z))\Big\vert_{z = z^{\ast}}
\\
& {}
  + c\, ( U(z) - U(z^{*}) )
  + \int_{z^{*}}^z g(U(z')) \,\mathrm{d}z' = 0 \,.
\end{aligned}
\end{equation}
\end{itemize}
\endgroup
\end{definition}
\par\vskip 10pt

\begin{remark}\label{rem-def_TW}\nopagebreak
\begingroup\rm
An important feature of our definition of
a travelling wave for problem \eqref{e:FKPP} above is the fact
that we do {\em not\/} assume that its profile,
$U\colon \RR\to \RR$, is a sufficiently smooth function that obeys
the differential equation \eqref{eq:FKPP} in a classical sense.
In fact, we will see in the next remark
(Remark~\ref{rem-trav_wave}, Part~{\rm (ii)})
that the ``weaker'' integral form of eq.~\eqref{eq:FKPP},
given in eq.~\eqref{eq:int_U(z)} above,
easily yields also the ``stronger'' classical form~\eqref{eq:FKPP}
at every point $z\in \RR$ such that $U(z)\not\in \{ 0,\, 1\}$.
In other words, in case the wave profile $U$ is only continuous,
but not differentiable, one has to take advantage of the integral form
\eqref{eq:int_U(z)} only for $z\in \RR$ near those points
${\tilde z}\in \RR$ at which $U({\tilde z})\in \{ 0,\, 1\}$.

The integral form \eqref{eq:int_U(z)} enables us to use
rather general, nonsmooth diffusion and reaction terms,
$d$ and $g$, respectively.
Last but not least, our definition of a travelling wave covers
both alternatives for travelling waves introduced in
{\sc F.\ S\'anchez\--Gardu\~{n}o} and {\sc P.~K.\ Maini}
\cite[Sect.~3, p.~167]{Sanchez-Maini}:
{\em front\--type\/} and {\em sharp\--type\/} travelling waves
(see also \cite[Sect.~2, pp.\ 473--474]{Malag-Marce-1}).
It has been shown in
{\sc L.\ Malaguti} and {\sc C.\ Marcelli}
\cite[Sect.~2, pp.\ 476--481]{Malag-Marce-1}
that the cases of
$z_0 > -\infty$ and/or $z_1 < +\infty$
may occur if the nonlinear reaction function $g(s)$ and
the diffusion term $d(s)$ are not differentiable at the points
$s\in \{ 0,\, 1\}$.
In accordance with
\cite[Remark~1, p.~478]{Malag-Marce-1},
we now persue the case of $g$ and/or $d$ being ``nonsmooth''
at the points $s\in \{ 0,\, 1\}$ in greater details.
\endgroup
\end{remark}
\par\vskip 10pt

Definition \ref{def-trav_wave} has the following simple,
but important technical consequences to be used in the sequel:

\begin{remark}\label{rem-trav_wave}\nopagebreak
\begingroup\rm
{\rm (i)}$\;$
Equation~\eqref{eq:FKPP} being translation invariant
($z\mapsto z + \zeta\colon \RR\to \RR$, for $\zeta\in \RR$ fixed),
we are allowed to choose the profile $U$ in such a way that
$U(0) = 1/2$.
This choice will determine the profile, $U$, uniquely if needed,
thanks to the strict monotonicity of the profile
throughout the open interval $(z_0,z_1)\subset \RR$, by $U'< 0$;
cf.\ Proposition~\ref{prop-monot_TW} below.
However, we do not assume $U(0) = 1/2$, in general, unless
we need the uniqueness of $U$ for a fixed speed $c\in \RR$.

{\rm (ii)}$\;$
{\rm Hypothesis\/} {\bf (H1)} combined with
Definition~\ref{def-trav_wave}, Part~{\rm (b)},
imply that, at every point
$\xi\in \RR$ with $U(\xi)\not\in \{ 0,\, 1\}$, we have
$d(U(\xi)) > 0$ and the derivative $U'(\xi)$ exists and satisfies
\begin{math}
  \frac{\mathrm{d}}{\mathrm{d}z}\, D(U(z)) \bigg\vert_{z = \xi}
  = d(U(\xi))\, U'(\xi) .
\end{math}

{\rm (iii)}$\;$
There exist two sequences
$\xi_n\in (n,n+1)$ and
$\xi_n^{*}\in (-n-1,-n)$; $n=1,2,3,\dots$, such that
\begin{equation}
\label{lim:U'(z)_infty}
  \frac{\mathrm{d}}{\mathrm{d}z}\, D(U(z))\bigg\vert_{z = \xi_n}
  \,\longrightarrow\, 0 \quad\mbox{ and }\quad
  \frac{\mathrm{d}}{\mathrm{d}z}\, D(U(z))\bigg\vert_{z = \xi_n^{*}}
  \,\longrightarrow\, 0 \quad\mbox{ as }\; n\to \infty \,.
\end{equation}
Indeed, we can apply the mean value theorem to
the (continuously differentiable) function
$z\mapsto D(U(z))\colon \RR\to \RR$ in each of the intervals
$(n,n+1)$ and $(-n-1,-n)$; $n=1,2,3,\dots$, to conclude that there are
$\xi_n\in (n,n+1)$ and
$\xi_n^{*}\in (-n-1,-n)$, such that
\begin{align*}
& D(U(n+1)) - D(U(n)) =
  \frac{\mathrm{d}}{\mathrm{d}z}\, D(U(z))\Big\vert_{z = \xi_n}
    \quad\mbox{ and }\quad
\\
& D(U(-n)) - D(U(-n-1)) =
  \frac{\mathrm{d}}{\mathrm{d}z}\, D(U(z))\Big\vert_{z = \xi_n^{*}} \,.
\end{align*}
The limits in eq.~\eqref{lim:U'(z)_infty} follow from
Definition~\ref{def-trav_wave} combined with
the limits in~\eqref{lim:FKPP}.
\endgroup
\end{remark}
\par\vskip 10pt

\section{Basic properties of a wave profile}
\label{s:Wave_Profile}

Throughout this section we assume that
$d,g\colon \RR\to \RR$ satisfy
{\rm Hypotheses\/} {\bf (H1)} and\/ {\bf (H2)}.
Next, we show that any wave profile $U\colon \RR\to \RR$
takes only values between $0$ and $1$.

\begin{lemma}\label{lemma-values_TW}
{\rm (Wave profile values.)}$\,$
Let\/
\begin{math}
  (x,t)\mapsto u(x,t) = U(x-ct)\colon \RR\times \RR_+\to \RR
\end{math}
be a TW with speed\/ $c\in \RR$ and profile $U\colon \RR\to \RR$.
Then we have $0\leq U(z)\leq 1$ for every\/ $z\in \RR$.
\end{lemma}

\par\vskip 10pt
\proof
We have $U(z)\geq 0$ for every $z\in \RR$, by
Definition~\ref{def-trav_wave}.
By contradiction to $U(z)\leq 1$ for every $z\in \RR$,
suppose there is a number $\xi\in \RR$ such that $U(\xi) > 1$.
We make use of the limits in \eqref{lim:FKPP}
to conclude that there are numbers $\xi_1, \xi_2\in \RR$ such that
$\xi_1 < \xi < \xi_2$ and
$U(\xi) > \min\{ U(\xi_1), U(\xi_2)\} > 1$.
We may choose $\xi_1$ and $\xi_2$, close enough to $\xi$,
in such a way that also
$U(z) > 1$ holds for every $z\in [\xi_1, \xi_2]$.
Denoting by $\xi_0\in [\xi_1, \xi_2]$
a (global) maximizer for the function $U$ over the compact interval
$[\xi_1, \xi_2]$, we arrive at $\xi_0\in (\xi_1, \xi_2)$,
$U(\xi_0)\geq U(\xi) > 1$, $U'(\xi_0) = 0$, and
\begin{equation*}
    d(U(z))\, U'(z) - d(U(\xi_0))\, U'(\xi_0)
  = {}- c\, ( U(z) - U(\xi_0) )
  - \int_{\xi_0}^z g(U(z')) \,\mathrm{d}z'
\end{equation*}
for all $z\in [\xi_1, \xi_2]$, by eq.~\eqref{eq:int_U(z)} and
Remark~\ref{rem-trav_wave}, Part~{\rm (i)}.
Since $U'(\xi_0) = 0$, the last equation entails
\begin{equation}
\label{eq:int_U(xi_0)}
    d(U(z))\, \frac{ U'(z) - U'(\xi_0) }{z - \xi_0}
  = {}- c\, \frac{ U(z) - U(\xi_0) }{z - \xi_0}
  - \frac{1}{z - \xi_0}
    \int_{\xi_0}^z g(U(z')) \,\mathrm{d}z'
\end{equation}
for all
$z\in [\xi_1, \xi_2]\setminus \{ \xi_0\}$.
We apply the mean value theorem to the right\--hand side of
eq.~\eqref{eq:int_U(xi_0)} to conclude that, for every
$z\in [\xi_1, \xi_2]$, $z\neq \xi_0$, there is a number
${\hat z}\in [\xi_1, \xi_2]$ between $\xi_0$ and $z$, such that
\begin{equation*}
    d(U(z))\, \frac{ U'(z) - U'(\xi_0) }{z - \xi_0}
  = {}- c\, U'({\hat z}) - g(U({\hat z})) \,.
\end{equation*}
Letting $z\to \xi_0$ we get also ${\hat z}\to \xi_0$
and, consequently, the second derivative
$U''(\xi_0)$ of $U$ at $\xi_0$ exists and satisfies
\begin{equation*}
    d(U(\xi_0))\, U''(\xi_0) = {}- c\, U'(\xi_0) - g(U(\xi_0))
  = {}- g(U(\xi_0)) > 0 \,,
\end{equation*}
where $d(U(\xi_0)) > 0$.
Hence, $U'(\xi_0) = 0$ and $U''(\xi_0) > 0$ show that
$\xi_0\in (\xi_1, \xi_2)$ is also
a strict local minimizer for the function $U$ in the open interval
$(\xi_1, \xi_2)$.
But this contradicts our construction of $\xi_0$
as a (global) maximizer for $U$ over $[\xi_1, \xi_2]$.

This proves $U(z)\leq 1$ for all $z\in \RR$.
\qed
\par\vskip 10pt

Now we are ready to calculate the wave speed $c$
explicitly from the wave profile $U$.

\begin{lemma}\label{lemma-speed_TW}
{\rm (Wave speed.)}$\,$
Let\/
\begin{math}
  (x,t)\mapsto u(x,t) = U(x-ct)\colon \RR\times \RR_+\to \RR
\end{math}
be a TW with speed\/ $c\in \RR$ and profile $U\colon \RR\to \RR$.
Then we have $0\leq U(z)\leq 1$ and\/
$g(U(z))\geq 0$ for every\/ $z\in \RR$, together with\/
\begin{equation}
\label{e:int_g(U(z))}
  0 < c = \int_{-\infty}^{+\infty} g(U(z')) \,\mathrm{d}z' < \infty \,.
\end{equation}
Moreover, eq.~\eqref{eq:int_U(z)} is equivalent with
\begin{equation}
\label{e:int_U(z)}
    \frac{\mathrm{d}}{\mathrm{d}z}\, D(U(z))
  + c\, U(z) - \int_{z}^{+\infty} g(U(z')) \,\mathrm{d}z' = 0
  \quad\mbox{ for all }\, z\in \RR \,.
\end{equation}
\end{lemma}

\par\vskip 10pt
\proof
We have $0\leq U(z)\leq 1$ for every $z\in \RR$,
by Lemma~\ref{lemma-values_TW},
which yields $g(U(z))$ $\geq 0$, by {\rm Hypothesis\/} {\bf (H2)}.

For every fixed $n=1,2,3,\dots$ we take the pair
$(z^{*},z) = (\xi_n^{*}, \xi_n)$
in eq.~\eqref{eq:int_U(z)}, where the latter pair has been specified
in Remark~\ref{rem-trav_wave}, Part~{\rm (ii)}.
Applying \eqref{lim:FKPP} and \eqref{lim:U'(z)_infty}
to eq.~\eqref{eq:int_U(z)} and letting $n\to \infty$, we arrive at
\begin{equation*}
  {}- c + \int_{-\infty}^{+\infty} g(U(z')) \,\mathrm{d}z' = 0 \,,
\end{equation*}
by the Lebesgue monotone convergence theorem.
This proves eq.~\eqref{e:int_g(U(z))} with $c\geq 0$.
However, the integrand $g(U(z'))\geq 0$ cannot vanish identically
for all $z'\in \RR$, by the continuity of the wave profile
$U\colon \RR\to \RR$ and the limits \eqref{lim:FKPP}
which guarantee $U({\hat z}) = \frac{1}{2}\in (0,1)$
for some ${\hat z}\in \RR$; hence, $g(U({\hat z})) > 0$.
Since also $g\colon \RR\to \RR$ is continuous, by
{\rm Hypothesis\/} {\bf (H2)}, we must have $c > 0$, by
eq.~\eqref{e:int_g(U(z))}.

To verify also eq.~\eqref{e:int_U(z)}, we now take the pair
$(z^{*},z) = (\xi_n^{*}, z)$, where
$\xi_n^{*}\in (-n-1,-n)$ is as above and $z\in \RR$ is arbitrary.
Applying \eqref{lim:FKPP} and \eqref{lim:U'(z)_infty}
to eq.~\eqref{eq:int_U(z)} again and letting $n\to \infty$, we obtain
\begin{equation*}
    \frac{\mathrm{d}}{\mathrm{d}z}\, D(U(z))
  + c\, ( U(z) - 1 )
  + \int_{-\infty}^{z} g(U(z')) \,\mathrm{d}z' = 0
  \quad\mbox{ for all }\, z\in \RR \,.
\end{equation*}
Finally, we apply \eqref{e:int_g(U(z))} to the last equation
to derive \eqref{e:int_U(z)}.
\qed
\par\vskip 10pt

We continue with the constant sections of the travelling wave.

\begin{lemma}\label{lemma-const_TW}
{\rm (Constant sections.)}$\,$
Let\/
\begin{math}
  (x,t)\mapsto u(x,t) = U(x-ct)\colon \RR\times \RR_+\to \RR
\end{math}
be a TW with speed\/ $c\in \RR$ and profile $U\colon \RR\to \RR$.
Assume that\/ $\xi\in \RR$ is such that\/ $U(\xi)\in \{ 0,\, 1\}$.
Then the following two alternatives are valid:
\begin{itemize}
\item[{\rm (i)}]
$\;$
If\/ $U(\xi) = 0$ then $U(z)\equiv 0$ for every\/ $z\geq \xi$.
\item[{\rm (ii)}]
$\;$
If\/ $U(\xi) = 1$ then $U(z)\equiv 1$ for every\/ $z\leq \xi$.
\end{itemize}
\end{lemma}

\par\vskip 10pt
\proof
We recall that $0\leq U(z)\leq 1$ for every $z\in \RR$,
by Lemma~\ref{lemma-values_TW}.

Alt.~(i):
$\;$
Assume that $U(\xi) = 0$.
Suppose there is some $\xi^{*}\in (\xi, +\infty)$ such that
$U(\xi^{*}) > 0$.
We can guarantee even $0 < U(\xi^{*}) < 1$,
by taking $\xi^{*}\in (\xi, +\infty)$ closer to $\xi$.
This implies $g(U(\xi^{*})) > 0$ and, consequently, we have
\begin{math}
  \int_{\xi}^{+\infty} g(U(z')) \,\mathrm{d}z' > 0 .
\end{math}
Furthermore, our definition of a travelling wave, 
Definition~\ref{def-trav_wave}, Part~{\rm (b)},
guarantees that also
\begin{math}
  \frac{\mathrm{d}}{\mathrm{d}z}\, D(U(z))\Big\vert_{z = \xi} = 0 ,
\end{math}
thanks to $U(\xi) = 0$.
We insert these facts into eq.~\eqref{e:int_U(z)}
with $z = \xi$, which yields
\begin{math}
  \int_{\xi}^{+\infty} g(U(z')) \,\mathrm{d}z' = 0 ,
\end{math}
a contradiction with the inequality (${} > 0$) above.

Alt.~(ii):
$\;$
Now assume $U(\xi) = 1$ and
suppose there is some $\xi^{*}\in (-\infty, \xi)$ such that
$U(\xi^{*}) < 1$.
Again, we can guarantee $0 < U(\xi^{*}) < 1$,
by taking $\xi^{*}\in (-\infty, \xi)$ closer to $\xi$.
This implies $g(U(\xi^{*})) > 0$ and, as above, we have
\begin{math}
  \int_{-\infty}^{\xi} g(U(z')) \,\mathrm{d}z' > 0 .
\end{math}
Definition~\ref{def-trav_wave}, Part~{\rm (b)},
guarantees also
\begin{math}
  \frac{\mathrm{d}}{\mathrm{d}z}\, D(U(z))\Big\vert_{z = \xi} = 0 ,
\end{math}
thanks to $U(\xi) = 1$.
We insert these facts into eq.~\eqref{e:int_U(z)}
with $z = \xi$, which yields
\begin{math}
  \int_{\xi}^{+\infty} g(U(z')) \,\mathrm{d}z' = c .
\end{math}
A comparison of this equality with eq.~\eqref{e:int_g(U(z))}
forces
\begin{math}
  \int_{-\infty}^{\xi} g(U(z')) \,\mathrm{d}z' = 0 ,
\end{math}
a contradiction with the inequality (${} > 0$) above.

The lemma is proved.
\qed
\par\vskip 10pt

Finally, we establish the monotonicity of the travelling wave
(see {\rm Definition~\ref{def-trav_wave}}).

\begin{proposition}\label{prop-monot_TW}
{\rm (Monotonicity.)}$\,$
Let\/
\begin{math}
  (x,t)\mapsto u(x,t) = U(x-ct)\colon \RR\times \RR_+\to \RR
\end{math}
be a TW with speed\/ $c\in \RR$ and profile $U\colon \RR\to \RR$.
Then $0\leq U(z)\leq 1$ and\/ $U'(z)\leq 0$ for all\/ $z\in \RR$.
Moreover, there is an open interval\/ $(z_0,z_1)\subset \RR$,
$-\infty\leq z_0 < z_1\leq +\infty$, such that\/
$U'< 0$ on $(z_0,z_1)$ together with
\begin{equation*}
\left\{\quad
  \begin{alignedat}{2}
& \lim_{z\to z_0+} U(z) = 1 \quad\mbox{ and }\quad
&&  U(z) = 1 \;\mbox{ if }\; -\infty < z\leq z_0 \,,
\\
& \lim_{z\to z_1-} U(z) = 0 \quad\mbox{ and }\quad
&&  U(z) = 0 \;\mbox{ if }\; z_1\leq z < +\infty \,.
  \end{alignedat}
\right.
\end{equation*}
\end{proposition}

\par\vskip 10pt
\proof
Recalling Lemmas
\ref{lemma-values_TW}, \ref{lemma-speed_TW}, and \ref{lemma-const_TW},
we conlude that it remains to prove
$U'(z) < 0$ for every $z\in \RR$ satisfying $0 < U(z) < 1$.
Suppose not; hence, there is some $\xi\in \RR$ such that
$U'(\xi) = 0$ and $0 < U(\xi) < 1$.
Eq.~\eqref{eq:int_U(z)} and
Remark~\ref{rem-trav_wave}, Part~{\rm (i)}, yield
\begin{equation}
\label{eq:int_v(xi_0)}
    d(U(z))\, U'(z) - d(U(\xi))\, U'(\xi)
  = {}- c\, ( U(z) - U(\xi) )
  - \int_{\xi}^z g(U(z')) \,\mathrm{d}z'
\end{equation}
for all $z\in \RR$, in analogy with our proof of
Lemma~\ref{lemma-values_TW}, eq.~\eqref{eq:int_U(xi_0)}.

Next, we show that every such point $\xi$ must be
a strict (i.e., isolated) local maximum satisfying $U''(\xi) < 0$.
Let us choose $\xi_1, \xi_2\in \RR$ such that
$\xi_1 < \xi < \xi_2$ and $0 < U(z) < 1$ holds for all
$z\in [\xi_1, \xi_2]$.
We apply the mean value theorem to the right\--hand side of
eq.~\eqref{eq:int_v(xi_0)} to conclude that, for every
$z\in [\xi_1, \xi_2]$, $z\neq \xi$, there is a number
${\hat z}\in [\xi_1, \xi_2]$ between $\xi$ and $z$, such that
\begin{equation*}
    d(U(z))\, \frac{ U'(z) - U'(\xi) }{z - \xi}
  = {}- c\, U'({\hat z}) - g(U({\hat z})) \,.
\end{equation*}
Letting $z\to \xi$ we conclude that ${\hat z}\to \xi$,
$d(U(z))\to d(U(\xi)) > 0$, and
\begin{equation*}
    d(U(\xi))\, U''(\xi) = {}- c\, U'(\xi) - g(U(\xi))
                         = {}- g(U(\xi)) < 0 \,.
\end{equation*}
This yields $U''(\xi) < 0$.

Since $U(z)\to 1$ as $z\to -\infty$, and $U(\xi) < 1$,
there is some $\xi_1'\in (-\infty, \xi)$ such that
$U(\xi) < U(\xi_1') < 1$.
Now let $\xi_0\in [\xi_1', \xi]$ be a (global) minimizer for
the function $U$ over the compact interval $[\xi_1', \xi]$.
With a help from $U'(\xi) = 0$ and $U''(\xi) < 0$,
we arrive at $\xi_0\in (\xi_1', \xi)$,
$U(\xi_0) < U(\xi) < 1$, $U'(\xi_0) = 0$, and
eq.~\eqref{eq:int_v(xi_0)} with $\xi_0$ in place of~$\xi$.
But then, by what we have proved above, if also $U(\xi_0) > 0$
then we must have $U''(\xi_0) < 0$ as above.
This contradicts our choice of $\xi_0$ to be a (global) minimizer for
the function $U$ over the open interval $(\xi_1', \xi)$.

The case $U(\xi_0) = 0$ would lead to a contradiction, by
Lemma~\ref{lemma-const_TW}.
It would force
$U(z) = 0$ for every $z\geq \xi_0$ and, in particular, also
$U(\xi) = 0$, thus contradicting our choice of $\xi\in \RR$.

We conclude that $U'(z) < 0$ holds for every $z\in (z_0,z_1)$.
\qed
\par\vskip 10pt

\section{A phase plane transformation}
\label{s:Phase_Plane}

We use a phase plane transformation
(cf.\ {\sc J.~D.\ Murray} \cite{Murray-I}, {\S}13.2, pp.\ 440--441,
 {\sc L.\ Malaguti} and {\sc C.\ Marcelli} \cite{Malag-Marce-1},
 {\sc R.\ Engui\c{c}a}, {\sc A.\ Gavioli}, and {\sc L.\ Sanchez}
 \cite[Sect.~1]{EGSanchez}, and
 {\sc P.\ Dr\'abek} and {\sc P.\ Tak\'a\v{c}} \cite{DrabTak-1})
in order to describe all monotone decreasing travelling waves
$u(x,t)\equiv U(x-ct - \zeta)$
where $U\colon \RR\to \RR$ is the profile of a travelling wave
normalized by $U(0) = 1/2$ as specified in
{\rm Remark~\ref{rem-trav_wave}}, {\rm Part~(i)}, and
$\zeta\in \RR$ is a suitable constant; see also
Proposition~\ref{prop-monot_TW}.
We reduce the second\--order differential equation for $U = U(z)$
to a first\--order ordinary differential equation for
the derivative
${\mathrm{d}z} / {\mathrm{d}U}$ of its inverse function
$U\mapsto z = z(U)$ as a function of $U\in (0,1)$.
In fact, below we find a nonlinear differential equation
for the derivative
\begin{equation*} 
  U'(z) = \genfrac{(}{)}{}0{\mathrm{d}z}{\mathrm{d}U}^{-1}
        \equiv \frac{1}{z'(U)} < 0
  \quad\mbox{ as a function of }\, U\in (0,1) \,.
\end{equation*}
To this end, we make the substitution
\begin{equation}
\label{e:dU/dz}
  V\eqdef {}- d(U)\, \genfrac{}{}{}0{\mathrm{d}U}{\mathrm{d}z} > 0
  \quad\mbox{ for }\, z\in (z_0,z_1)
\end{equation}
and consequently look for $V = V(U)$ as a function of
$U\in (0,1)$ that satisfies the following differential equation
obtained from eq.~\eqref{eq:FKPP}:
\begin{equation*}
{}- \frac{\mathrm{d}V}{\mathrm{d}U}\cdot \frac{\mathrm{d}U}{\mathrm{d}z}
  + c\, \frac{\mathrm{d}U}{\mathrm{d}z} + g(U)
  = 0 \,,
    \quad z\in (z_0,z_1) \,,
\end{equation*}
that is,
\begin{equation}
\label{eq:FKPP:V(U)}
    \frac{\mathrm{d}V}{\mathrm{d}U}\cdot \frac{V}{d(U)}
  - c\, \frac{V}{d(U)} + g(U) = 0 \,,
    \quad U\in (0,1) \,.
\end{equation}
Hence, we are looking for the inverse function
$U\mapsto z(U)$ with the derivative
\begin{equation*}
    \genfrac{}{}{}0{\mathrm{d}z}{\mathrm{d}U}
  = {}- \frac{d(U)}{V(U)} < 0
  \quad\mbox{ for }\, U\in (0,1) \,, \quad\mbox{ such that }\,
  z(1/2) = 0 \,.
\end{equation*}
Finally, we multiply eq.~\eqref{eq:FKPP:V(U)} by $d(U)$,
make the substitution
\begin{equation}
\label{e:y=V^p'}
  y = V^2
  = d(U)^2\, \genfrac{|}{|}{}0{\mathrm{d}U}{\mathrm{d}z}^2
  = \left| \frac{\mathrm{d}}{\mathrm{d}z}\, D(U(z))
    \right|^2 > 0 \,,
\end{equation}
and write $r$ in place of $U$, thus arriving at
\begin{equation*}
    \frac{1}{2}\cdot \frac{\mathrm{d}y}{\mathrm{d}r}
  - c\, \sqrt{y} + f(r) = 0 \,,
    \quad r\in (0,1) \,.
\end{equation*}
Here, the function $f\colon \RR\setminus \{0,\, 1\}\to \RR$
is defined by $f(r)\eqdef d(r)\, g(r)$ for every
$r\in \RR\setminus \{0,\, 1\}$.
Observe that $f$ is continuous on $\RR\setminus \{0,\, 1\}$ with
$f(r) > 0$ for every $r\in (0,1)$, and
$f(r) < 0$ for every $r\in (-\infty,0)\cup (1,\infty)$.
In our existence results in {\S}\ref{ss:exist_TW}
we will assume also
$\lim_{r\to 0+} f(r) = 0$ and $\lim_{r\to 1-} f(r) = 0$, that is,
the restriction $f|_{(0,1)}$ of $f$ to the open interval $(0,1)$
can be extended to a continuous function $f|_{[0,1]}$ on $[0,1]$
by setting $f(0) = f(1) = 0$.

This means that the unknown function
$y\colon (0,1)\to (0,\infty)$ of $r$ verifies also
\begin{equation}
\label{eq:FKPP:y(r)}
    \frac{\mathrm{d}y}{\mathrm{d}r}
  = 2\left( c\, \sqrt{y^{+}} - f(r) \right) \,,
    \quad r\in (0,1) \,,
\end{equation}
where $y^{+} = \max\{ y,\, 0\}$.
Since we require that the function
$z\mapsto D(U(z))\colon \RR\to \RR$
be continuously differentiable with the derivative
\begin{math}
  \frac{\mathrm{d}}{\mathrm{d}z}\, D(U(z))
\end{math}
vanishing at every point\/ $\xi\in \RR$ such that\/
$U(\xi)\in \{ 0,\, 1\}$, that is,
$\frac{\mathrm{d}}{\mathrm{d}z}\, D(U(z)\to 0$
as $z\to z_0+$ and $z\to z_1-$, the function
$y = y(r) = | \mathrm{d} D(U(z)) / \mathrm{d}z |^2$
must satisfy the boundary conditions
\begin{equation}
\label{bc:FKPP:y(r)}
  y(0) = y(1) = 0 \,.
\end{equation}

The results of our phase plane transformation are collected
in the following lemma.
Recall that, by
Lemma~\ref{lemma-speed_TW}, eq.~\eqref{e:int_g(U(z))},
any TW with speed $c\in \RR$, if it exists, must have speed $c > 0$.

\begin{lemma}\label{lem-y_c}
{\rm (Existence of the wave profile.)}$\,$
Assume that\/ $d$ and\/ $g$ satisfy\/
{\rm Hypotheses\/} {\bf (H1)} and\/ {\bf (H2)}, respectively.
Let\/ $c\in (0,\infty)$.
Then problem \eqref{eq:FKPP:y(r)}, \eqref{bc:FKPP:y(r)}
has a classical solution
\begin{math}
  y\equiv y_{\mathrm{c}}\colon (0,1)\to (0,\infty)
\end{math}
if and only if problem \eqref{eq:FKPP}, \eqref{z_0,1:FKPP}
has a solution
\begin{math}
  U\colon (z_0,z_1)\to (0,\infty) \,.
\end{math}
\end{lemma}
\par\vskip 10pt

In the next two paragraphs we are concerned with
the solvability of the overdetermined first\--order
boundary value problem
\eqref{eq:FKPP:y(r)}, \eqref{bc:FKPP:y(r)}
with a free parameter $c\in \RR$.
We address the natural questions, such as
existence and non\-existence, and
uniqueness and non\-uniqueness of a classical solution
$y\colon (0,1)\to (0,\infty)$.

\subsection{A nonexistence result}
\label{ss:nonexist_TW}

In contrast with the well\--known existence results due to
{\sc L.\ Malaguti} and {\sc C.\ Marcelli} \cite{Malag-Marce-1} and
{\sc R.\ Engui\c{c}a}, {\sc A.\ Gavioli}, and {\sc L.\ Sanchez}
 \cite{EGSanchez},
we now formulate a nonexistence result for a TW
$u(x,t)\equiv U(x-ct - \zeta)$ whose profile
$U\colon \RR\to \RR$ should satisfy the boundary value problem
\eqref{eq:FKPP}, \eqref{z_0,1:FKPP}.

\begin{proposition}\label{prop-nonexist_TW}
{\rm (Nonexistence of TW.)}$\,$
Let speed\/ $c\in (0,\infty)$ be arbitrary and assume that\/
$f = dg\colon \RR\setminus \{0,\, 1\}\to \RR$
satisfies the following {\em growth rate\/} condition,
\begin{equation}
\label{ineq:f(r)/r,r=0}
  {f(r)} / {r}\geq c^2 \quad\mbox{ for all }\, r\in (0,\delta) \,,
\end{equation}
where $\delta\in (0,1)$ is some number.
Then problem
\eqref{eq:FKPP:y(r)}, \eqref{bc:FKPP:y(r)}
has no classical solution $y\colon (0,1)\to \RR$.

In particular, if $f$ satisfies
\begin{equation}
\label{lim:f(r)/r,r=0}
  \lim_{r\to 0+} \,\frac{f(r)}{r} = +\infty \,,
\end{equation}
then problem
\eqref{eq:FKPP:y(r)}, \eqref{bc:FKPP:y(r)}
has no classical solution $y\colon (0,1)\to \RR$ for any\/ $c\in \RR$.
Consequently, also problem \eqref{e:FKPP}
has no TW solution in the sense of\/
{\rm Definition~\ref{def-trav_wave}}.
\end{proposition}

\par\vskip 10pt
\proof
On the contrary, assume that
$y\colon (0,1)\to \RR$ is a classical solution to problem
\eqref{eq:FKPP:y(r)}, \eqref{bc:FKPP:y(r)}.
Then it must satisfy
$y(r) > 0$ for every $r\in (0,1)$.
Indeed, suppose that $y(r_0)\leq 0$ for some $r_0\in (0,1)$.
Owing to the zero boundary conditions \eqref{bc:FKPP:y(r)},
we may assume that $y$ attains its global minimum at $r_0$, i.e.,
$y(r_0) = \min_{r\in (0,1)} y(r)$.
Hence, we get $y'(r_0) = 0$.
But then eq.~\eqref{eq:FKPP:y(r)} at $r = r_0$ forces
\begin{equation*}
  0 = y'(r_0) - 2c\, \sqrt{y^{+}(r_0)}
    = {}- 2\, f(r_0) < 0 \,,
\end{equation*}
a contradiction.
We conclude that eq.~\eqref{eq:FKPP:y(r)} is equivalent with
\begin{equation}
\label{eq:FKPP:y(r)>0}
    \frac{\mathrm{d}}{\mathrm{d}r}\, \sqrt{y(r)}
  = c - \frac{f(r)}{ \sqrt{y(r)} } \,,
    \quad r\in (0,1) \,.
\end{equation}

Given $s\in (0,1)$ arbitrary, we integrate the inequality
\begin{math}
  \frac{\mathrm{d}}{\mathrm{d}r}\, \sqrt{y(r)}\leq c
\end{math}
over the interval $[0,s]$, thus arriving at
$\sqrt{y(s)}\leq \sqrt{y(0)} + cs = cs$ or, equivalently,
$y(r)\leq c^2 r^2$ for all $r\in (0,1)$.
We apply this estimate to the right\--hand side of
eq.~\eqref{eq:FKPP:y(r)>0} to conclude that
\begin{equation}
\label{ineq:FKPP:y(r)>0}
    \frac{\mathrm{d}}{\mathrm{d}r}\, \sqrt{y(r)}
  \leq c - \frac{f(r)}{cr}
  = \frac{1}{c} \left( c^2 - \frac{f(r)}{r} \right) \,,
    \quad r\in (0,1) \,.
\end{equation}
Finally, we apply ineq.~\eqref{ineq:f(r)/r,r=0}
to the right\--hand side of
ineq.~\eqref{ineq:FKPP:y(r)>0} for $r\in (0,\delta)$ to get
\begin{math}
  \frac{\mathrm{d}}{\mathrm{d}r}\, \sqrt{y(r)}\leq 0
\end{math}
for every $r\in (0,\delta)$.
This yields
$\sqrt{y(r)}\leq \sqrt{y(0)} = 0$ for every $r\in (0,\delta)$,
thus forcing $y(r)\equiv 0$ for $r\in [0,\delta]$,
a contradiction to $y(r) > 0$ for all $r\in (0,1)$.

We have proved that problem
\eqref{eq:FKPP:y(r)}, \eqref{bc:FKPP:y(r)}
has no classical solution $y\colon [0,1]\to \RR$.
The remaining statements of the proposition follow easily.
\qed
\par\vskip 10pt

\subsection{Some existence results}
\label{ss:exist_TW}

The following existence result for problem
\eqref{eq:FKPP:y(r)}, \eqref{bc:FKPP:y(r)} is due to
{\sc R.\ Engui\c{c}a}, {\sc A.\ Gavioli}, and {\sc L.\ Sanchez}
\cite[Prop.~2, p.~176]{EGSanchez}.

\begin{proposition}\label{prop-exist_TW}
{\rm (Existence of TW.)}$\,$
Assume that\/
$f = dg\colon \RR\setminus \{0,\, 1\}\to \RR$
satisfies $f(1) = 0$ and\/
\begin{equation}
\label{sup:f(r)/r,r=0}
  0 < \mu\eqdef \sup_{r\in (0,1)} \,\frac{f(r)}{r} < +\infty \,.
\end{equation}
Then there exists a number\/
$c^{*}\in \bigl( 0,\, 2\, \sqrt{\mu} \bigr]$ such that problem\/
\eqref{eq:FKPP:y(r)}, \eqref{bc:FKPP:y(r)} with speed $c\in \RR$
admits a unique positive solution if and only if\/
$c\geq c^{*}$.
Consequently, also problem \eqref{e:FKPP}
has a TW solution in the sense of\/
{\rm Definition~\ref{def-trav_wave}}.
\end{proposition}
\par\vskip 10pt

Some more related existence results can be found in
{\sc L.\ Malaguti} and {\sc C.\ Marcelli}
\cite[Theorems 2 and~3, pp.\ 474--475]{Malag-Marce-1}
for travelling waves distinguished by
the {\em front-\/} or {\em sharp\--type\/};
see our Figures 1, 2 or~3, respectively.

\begin{figure}[h]
\begin{center}
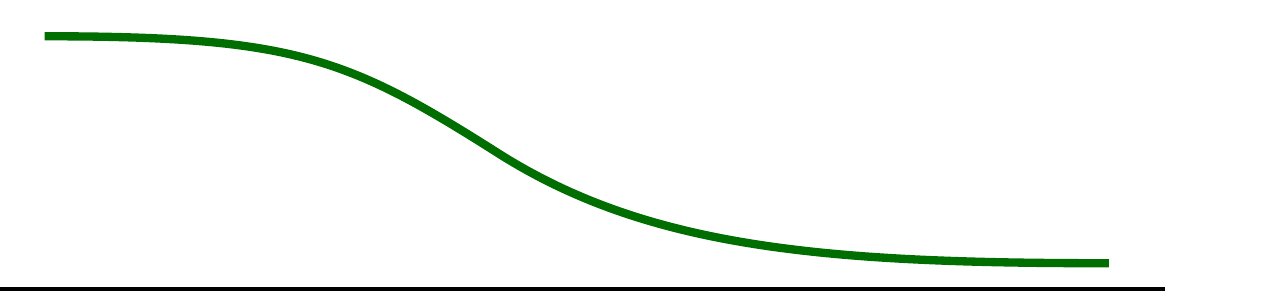
\caption{Travelling wave of front-type with $z_0=-\infty$, $z_1=+\infty$.}
\end{center}
\end{figure}

\begin{figure}[h]
\begin{center}
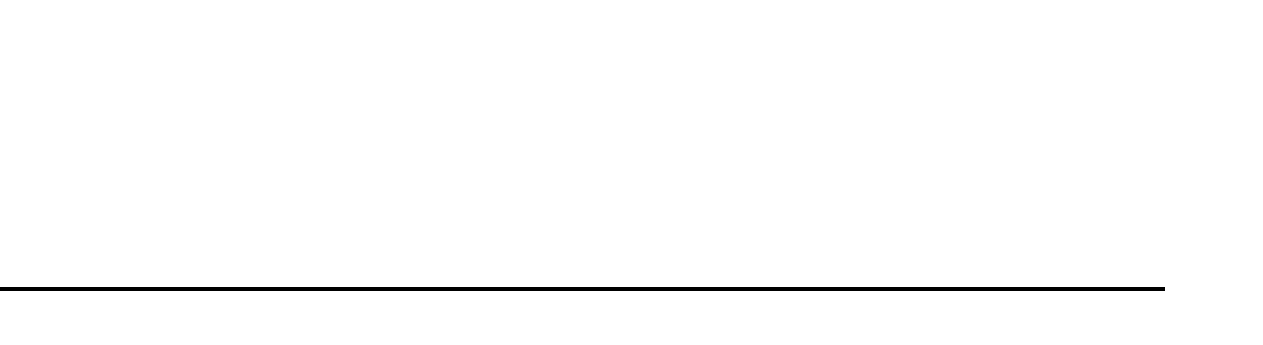
\caption{Travelling wave of front-type with $z_0>-\infty$, $z_1=+\infty$.}
\end{center}
\end{figure}

\begin{figure}[h]
\begin{center}
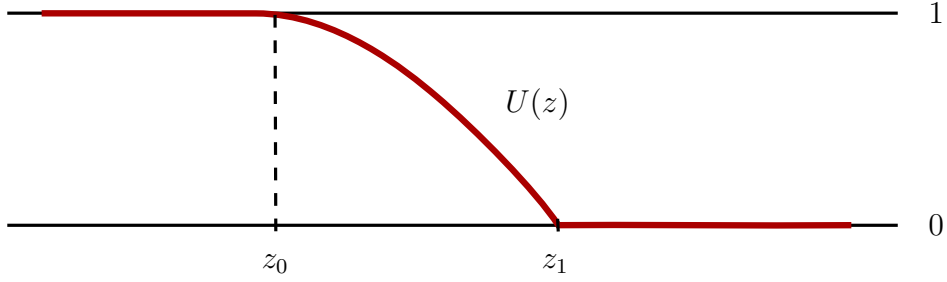
\caption{Travelling wave of sharp-type with $z_0>-\infty$, $z_1<+\infty$.}
\end{center}
\end{figure}

\begin{remark}\label{rem-rel:d,g}\nopagebreak
\begingroup\rm
Notice that conditions
\eqref{lim:f(r)/r,r=0} and \eqref{sup:f(r)/r,r=0}
impose a restriction on the mutual relation between
the diffusion $d(r)$ and the reaction $g(r)$ as $r\to 0+$.
In particular, given a reaction function $g\colon \RR\to \RR$
satisfying {\rm Hypothesis\/} {\bf (H2)},
diffusion $d(r)$ that degenerates ``suitably'' as $r\to 0+$
may guarantee the {\em existence\/} of a solution to problem
\eqref{eq:FKPP:y(r)}, \eqref{bc:FKPP:y(r)}.
On the other hand, diffusion $d(r)$ that blows up ``suitably''
as $r\to 0+$
may prevent the {\em existence\/} of a solution to
\eqref{eq:FKPP:y(r)}, \eqref{bc:FKPP:y(r)}.
\endgroup
\end{remark}
\par\vskip 10pt

\section{Interaction between diffusion and reaction,\\
         asymptotic shape of travelling waves}
\label{s:Trav-Shape}

In this section we prove a number of specialized results on
the profile of a travelling wave for some simple forms of
the nonlinearities $d(r)$ and $g(r)$ involved.
Our main goal here is to illustrate the biological meaning of
our mathematical results rather than to treat
mathematically general cases.
We restrict ourselves to diffusion and reaction terms
$d(r)$ and $g(r)$ having the following
{\it power\--type asymptotic behavior\/}
as $r\to 0+$ and $r\to 1-$, respectively, where
$\gamma_0$, $\gamma_1$, $\delta_0$, and $\delta_1$
are some real constants:
\begin{equation}
\label{e:d(1-)}
\left\{\qquad
\begin{alignedat}{2}
  \lim_{r\to 0+} \frac{g(r)}{r^{\gamma_0}}
& {}\eqdef g_0\in (0,\infty) \,,
\\
  \lim_{r\to 1-} \frac{g(r)}{(1-r)^{\gamma_1}}
& {}\eqdef g_1\in (0,\infty) \,,
\\
  \lim_{r\to 0+} \frac{d(r)}{r^{\delta_0}}
& {}\eqdef d_0\in (0,\infty) \,,
\\
  \lim_{r\to 1-} \frac{d(r)}{(1-r)^{\delta_1}}
& {}\eqdef d_1\in (0,\infty) \,.
\end{alignedat}
\right.
\end{equation}

The following restrictions on the parameters
$\gamma_0$, $\gamma_1$, $\delta_0$, and $\delta_1$
are imposed by {\rm Hypotheses\/} {\bf (H1)} and\/ {\bf (H2)}:
\begin{alignat*}{2}
& \mbox{ {\rm Hypothesis\/} {\bf (H1)} }\quad \Longrightarrow\quad
&&\delta_0 > -1 \;\mbox{ and }\; \delta_1 > -1 \,,
\\
& \mbox{ {\rm Hypothesis\/} {\bf (H2)} }\quad \Longrightarrow\quad
&&\gamma_0 > 0 \;\mbox{ and }\; \gamma_1 > 0 \,.
\end{alignat*}
In addition, recalling $f(r) = d(r)\, g(r)$ for every
$r\in \RR\setminus \{0,\, 1\}$, and $f$ continuous on $[0,1]$
with $f(0) = f(1) = 0$, we get also the restrictions
\begin{equation*}
  \gamma_0 + \delta_0 > 0 \quad\mbox{ and }\quad
  \gamma_1 + \delta_1 > 0 \,.
\end{equation*}

In what follows we treat the profile of the travelling wave
$r = U(z)$ for values near the equilibrium points
$r=0$ (in {\S}\ref{ss:asympt_0}) and
$r=1$ (in {\S}\ref{ss:asympt_1}).

\subsection{Profile asymptotics near $0$.}
\label{ss:asympt_0}

Let us define the following parameter sets,
see Figure~4,
\begin{alignat*}{2}
& \mathscr{M}_0^1\eqdef
  \{ (\gamma_0, \delta_0)\in \RR^2\colon
     \gamma_0 > 0 ,\ \delta_0 > -1 ,\ 0 < \gamma_0 + \delta_0 < 1
  \} \,,
\\
& \mathscr{M}_0^2\eqdef
  \{ (\gamma_0, \delta_0)\in \RR^2\colon
     \gamma_0 > 0 ,\ \delta_0 > -1 ,\ \gamma_0 + \delta_0\geq 1
  \} \,.
\end{alignat*}

\begin{figure}[h]
\begin{center}
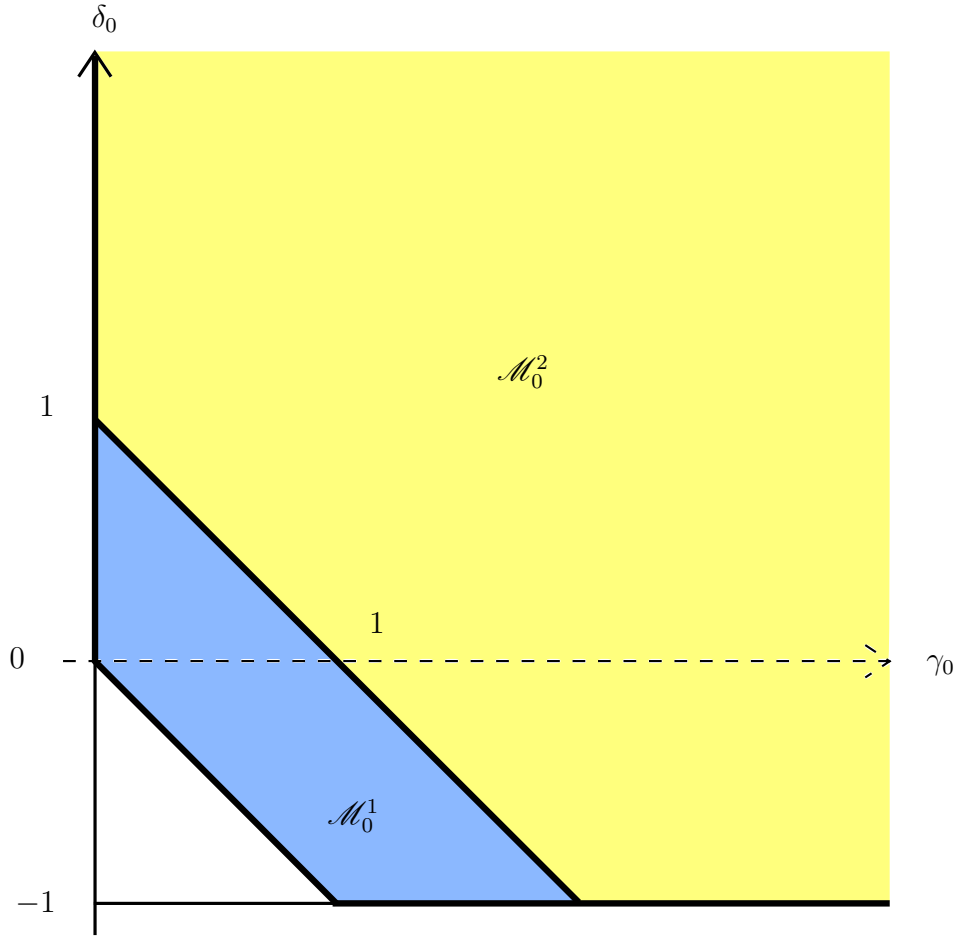
\caption{The sets $\mathscr{M}_0^1$ and $\mathscr{M}_0^2$.}
\end{center}
\end{figure}

For the parameter pairs
$(\gamma_0, \delta_0)\in \mathscr{M}_0^1\cup \mathscr{M}_0^2$
we have the following conclusions on the existence of travelling waves;
see {\rm Propositions\/} \ref{prop-nonexist_TW} and~\ref{prop-exist_TW}
above for further details.

\begin{theorem}\label{thm-asympt_0}
{\rm (i)}$\;$
$(\gamma_0, \delta_0)\in \mathscr{M}_0^1$
implies eq.~\eqref{lim:f(r)/r,r=0} and, hence,
no travelling wave exists, by
{\rm Proposition~\ref{prop-nonexist_TW}}.

{\rm (ii)}$\;$
$(\gamma_0, \delta_0)\in \mathscr{M}_0^2$
implies eq.~\eqref{sup:f(r)/r,r=0} and, hence,
a travelling wave exists, by
{\rm Proposition~\ref{prop-exist_TW}}.
\end{theorem}
\par\vskip 10pt

\subsection{Profile asymptotics near $1$.}
\label{ss:asympt_1}

Here, we need the following parameter sets,
see Figure~5,
\begin{alignat*}{2}
& \mathscr{M}_1^1\eqdef
  \{ (\gamma_1, \delta_1)\in \RR^2\colon
     0 < \gamma_1 < 1 + \delta_1 ,\ 0 < \gamma_1 + \delta_1\leq 1
  \} \,,
\\
& \mathscr{M}_1^2\eqdef
  \{ (\gamma_1, \delta_1)\in \RR^2\colon
     0 < 1 + \delta_1\leq \gamma_1 ,\ 0 < \gamma_1 + \delta_1\leq 1
  \} \,,
\\
& \mathscr{M}_1^3\eqdef
  \{ (\gamma_1, \delta_1)\in \RR^2\colon
     0 < \gamma_1 < 1 ,\ \gamma_1 + \delta_1 > 1
  \} \,,
\\
& \mathscr{M}_1^4\eqdef
  \{ (\gamma_1, \delta_1)\in \RR^2\colon
     \gamma_1\geq 1 ,\ \delta_1 > -1 ,\ \gamma_1 + \delta_1 > 1
  \} \,.
\end{alignat*}

\begin{figure}[h]
\begin{center}
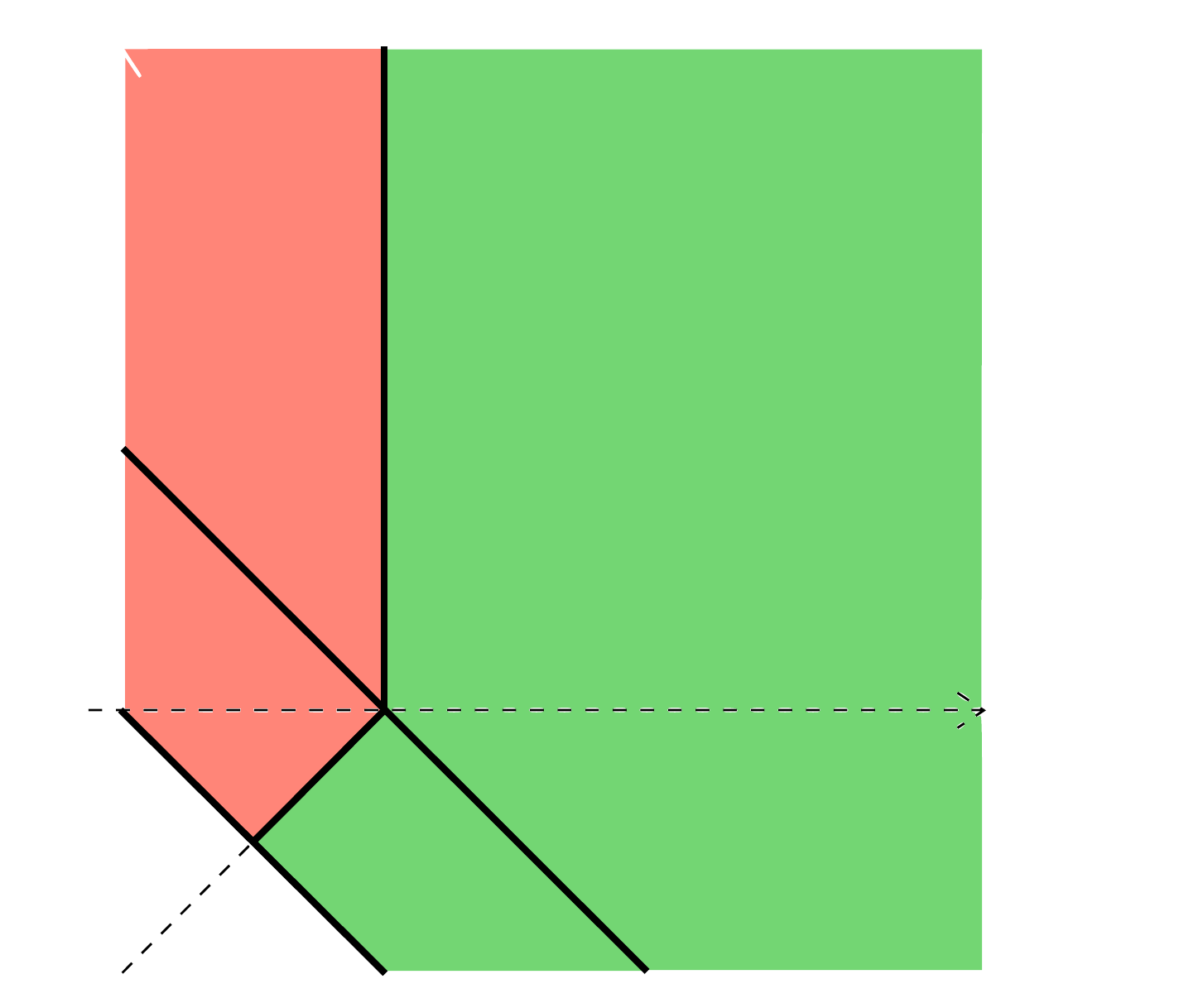
\caption{The sets $\mathscr{M}_1^1$, $\mathscr{M}_1^2$,
         $\mathscr{M}_1^3$, and $\mathscr{M}_1^4$.}
\end{center}
\end{figure}

From Section~\ref{s:Phase_Plane} we recall that
\begin{math}
  r\mapsto y\equiv y_{\mathrm{c}}(r)\colon (0,1)\to (0,\infty)
\end{math}
is a classical solution of problem
\eqref{eq:FKPP:y(r)}, \eqref{bc:FKPP:y(r)}.

In what follows we assume
$(\gamma_0, \delta_0)\in \mathscr{M}_0^2$, i.e.,
$y\equiv y_{\mathrm{c}}(r)$ exists
as a solution to the non\-linear two\--point boundary value problem
\eqref{eq:FKPP:y(r)}, \eqref{bc:FKPP:y(r)} for the unknown function
$y\colon (0,1)\to (0,\infty)$ with some speed $c > 0$, by
{\rm Theorem~\ref{thm-asympt_0}}{\rm (ii)} and
{\rm Proposition~\ref{prop-exist_TW}}.
Consequently, a travelling wave with the profile
$U\colon z\mapsto U(z)$ is obtained by the phase plane transformation
described in Section~\ref{s:Phase_Plane}, Lemma~\ref{lem-y_c}.
We classify the parameters $\gamma_1$ and $\delta_1$ according to
whether $z_0 > -\infty$ or $z_0 = -\infty$.
For the parameter pairs
$(\gamma_1, \delta_1)\in \cup_{i=1}^4 \mathscr{M}_1^i$
we will prove the following conclusions.

\begin{theorem}\label{thm-asympt_1}
Assume
$(\gamma_0, \delta_0)\in \mathscr{M}_0^2$.
Then we have:

{\rm (i)}$\;$
$z_0 > -\infty$ provided\/
$(\gamma_1, \delta_1)\in \mathscr{M}_1^1\cup \mathscr{M}_1^3$.

{\rm (ii)}$\;$
$z_0 = -\infty$ provided\/
$(\gamma_1, \delta_1)\in \mathscr{M}_1^2\cup \mathscr{M}_1^4$.
\end{theorem}

\par\vskip 10pt
\proof
We begin with

{\it Case 1:\/}$\quad$
$(\gamma_1, \delta_1)\in \mathscr{M}_1^1$.
We will compare the classical solution
$y\equiv y_{\mathrm{c}}\colon (0,1)\to (0,\infty)$
specified above with the function
\begin{math}
  w_{\kappa}(r)\eqdef \kappa\, (1-r)^{\gamma_1 + \delta_1 + 1}
\end{math}
of $r\in [0,1]$, where $\kappa > 0$
is a suitable number to be determined later.
We set $f_1 = d_1 g_1$ (${} > 0$) and write
\begin{math}
  f(r) = \left( f_1 + \eta(r)\right) (1-r)^{\gamma_1 + \delta_1} ,
\end{math}
where $\eta\colon [0,1]\to \RR$ is a continuous function with
$\eta(1) = 0$.

Then the differential operator in eq.~\eqref{eq:FKPP:y(r)}
takes the form
\begin{equation}
\label{oper:FKPP:y(r)}
  \mathscr{A}(y)(r)\eqdef
    \frac{\mathrm{d}y}{\mathrm{d}r} - 2 c\, \sqrt{y^{+}} + 2 f(r) \,,
    \quad r\in (0,1) \,.
\end{equation}
In particular, for the function $w_{\underline{\kappa}}$
defined above, with $\underline{\kappa} > 0$ small enough, we calculate
\begin{equation*}
\begin{aligned}
  \mathscr{A}(w_{\underline{\kappa}})(r)
  =
& {}- \underline{\kappa}
      (\gamma_1 + \delta_1 + 1)\, (1-r)^{\gamma_1 + \delta_1}
  - 2 c\, \sqrt{\underline{\kappa}}\,
      (1-r)^{ (\gamma_1 + \delta_1 + 1) / 2 }
\\
& {}+ 2 \left( f_1 + \eta(r)\right) (1-r)^{\gamma_1 + \delta_1} \,,
    \quad r\in (0,1) \,.
\end{aligned}
\end{equation*}
Since
$(\gamma_1, \delta_1)\in \mathscr{M}_1^1$
implies
\begin{math}
  \gamma_1 + \delta_1\leq \frac{1}{2} (\gamma_1 + \delta_1 + 1) ,
\end{math}
the first and third terms above dominate the second one
in the following sense, for $r\in (0,1)$ close enough to~$1$:
\begin{equation}
\label{oper:FKPP:w(r)}
\begin{aligned}
& \mathscr{A}(w_{\underline{\kappa}})(r)
  = (1-r)^{\gamma_1 + \delta_1}
\\
& \times
    \left[ {}- \underline{\kappa} (\gamma_1 + \delta_1 + 1)
             + 2 \left( f_1 + \eta(r)\right)
             - 2 c\, \sqrt{\underline{\kappa}}\,
               (1-r)^{ (1 - \gamma_1 - \delta_1) / 2 }
    \right] \,,
\end{aligned}
\end{equation}
provided $\underline{\kappa} > 0$ is chosen small enough,
relative to $f_1 > 0$.
This way we are able to guarantee
\begin{equation*}
  \mathscr{A}(w_{\underline{\kappa}})(r)
  \geq f_1\, (1-r)^{\gamma_1 + \delta_1} > 0
  \quad\mbox{ for all\/ $r\in (0,1)$ close to~$1$. }
\end{equation*}
Hence, there is a sufficiently small number $\varrho\in (0,1)$ such that
\begin{math}
  w_{\underline{\kappa}}\colon r\mapsto w_{\underline{\kappa}}(r)
\end{math}
is a {\it subsolution\/}
for the backward initial value problem
\begin{equation}
\label{back:FKPP:y(r)}
    \frac{\mathrm{d}y}{\mathrm{d}r}
  = 2\left( c\, \sqrt{y^{+}} - f(r) \right) \,,
    \quad r\in (1 - \varrho, 1) \,;\qquad y(1) = 0 \,.
\end{equation}

Recall that $c > 0$.
Observing that the nonlinearity $y\mapsto \sqrt{y^{+}}$
is a monotone, noninreasing function, we conclude that
the backward initial value problem \eqref{back:FKPP:y(r)}
possesses a unique classical solution
$y\equiv y_{\mathrm{c}}(r)$ on the interval $(1 - \varrho, 1)$.
By a similar monotonicity argument, we arrive at
\begin{math}
  y_{\mathrm{c}}(r)\geq w_{\underline{\kappa}}(r)
  = \underline{\kappa}\, (1-r)^{\gamma_1 + \delta_1 + 1}
\end{math}
for all $r\in (1 - \varrho, 1)$.
After returning to the original variables from eqs.\
\eqref{e:dU/dz} and \eqref{e:y=V^p'} we obtain
\begin{equation*}
  V(U)\geq \sqrt{\underline{\kappa}}\,
           (1-U)^{ (\gamma_1 + \delta_1 + 1) / 2 }
  \quad\mbox{ for all }\, U\in (1 - \varrho, 1) \,.
\end{equation*}
We combine this inequality with the last limit in \eqref{e:d(1-)}
to conclude that there is a constant $c_1 > 0$ such that
\begin{equation}
\label{ineq:dU/dz}
{}- \genfrac{}{}{}0{\mathrm{d}z}{\mathrm{d}U}
  = \frac{d(U)}{V(U)} \leq
    \frac{c_1}{ (1-U)^{ (\gamma_1 - \delta_1 + 1) / 2 } }
  \quad\mbox{ for all }\, U\in (1 - \varrho, 1) \,.
\end{equation}
Notice that the relation
$(\gamma_1, \delta_1)\in \mathscr{M}_1^1$
implies also
\begin{math}
  \frac{1}{2} (\gamma_1 - \delta_1 + 1) < 1 .
\end{math}
We fix an arbitrary number ${\tilde U}\in (1 - \varrho, 1)$, denote
${\tilde z} = z({\tilde U})\in (z_0,z_1)$ with
$U\mapsto z(U)\colon (0,1)\to (z_0,z_1)$
being the inverse function of $U\colon (z_0,z_1)\to (0,1)$,
and integrate ineq.~\eqref{ineq:dU/dz}
with respect to $U\in ({\tilde U}, 1)$, thus arriving at
\begin{equation*}
    {\tilde z} - z_0 = \int_{z_0}^{\tilde z} \mathrm{d}z
  = \int_1^{\tilde U}
    \genfrac{}{}{}0{\mathrm{d}z}{\mathrm{d}U} \,\mathrm{d}U
  = - \int_{\tilde U}^1
    \genfrac{}{}{}0{\mathrm{d}z}{\mathrm{d}U} \,\mathrm{d}U
  \leq c_1\int_{\tilde U}^1
    \frac{ \mathrm{d}U }{ (1-U)^{ (\gamma_1 - \delta_1 + 1) / 2 } }
  < \infty \,.
\end{equation*}
This estimate forces $z_0 > -\infty$.

{\it Case 2:\/}$\quad$
$(\gamma_1, \delta_1)\in \mathscr{M}_1^3$.
Here we compare
$y\equiv y_{\mathrm{c}}\colon (0,1)\to (0,\infty)$
with the new function
\begin{math}
  w_{\kappa}(r)\eqdef \kappa\, (1-r)^{ 2 (\gamma_1 + \delta_1) }
\end{math}
of $r\in [0,1]$, where $\kappa > 0$
is a suitable number to be determined later again.
Using eq.~\eqref{oper:FKPP:y(r)}, for
$\underline{\kappa} > 0$ small enough, we calculate
\begin{equation}
\label{ope:FKPP:w(r)}
\begin{aligned}
  \mathscr{A}(w_{\underline{\kappa}})(r)
  =
& {}- 2\underline{\kappa} (\gamma_1 + \delta_1)\,
  (1-r)^{ 2 (\gamma_1 + \delta_1) - 1 }
  - 2 c\, \sqrt{\underline{\kappa}}\, (1-r)^{ \gamma_1 + \delta_1 }
\\
& {}+ 2 \left( f_1 + \eta(r)\right) (1-r)^{\gamma_1 + \delta_1} \,,
    \quad r\in (0,1) \,.
\end{aligned}
\end{equation}
Since
$(\gamma_1, \delta_1)\in \mathscr{M}_1^3$
implies
\begin{math}
  2 (\gamma_1 + \delta_1) - 1 > \gamma_1 + \delta_1 ,
\end{math}
the second and third terms above dominate the first one
in the following sense, for $r\in (0,1)$ close enough to~$1$:
\begin{equation*}
  \mathscr{A}(w_{\underline{\kappa}})(r)
  \geq f_1\, (1-r)^{\gamma_1 + \delta_1} > 0 \,,
\end{equation*}
provided $\underline{\kappa} > 0$ is chosen small enough,
relative to $f_1 > 0$.
Hence, there is a sufficiently small number $\varrho\in (0,1)$ such that
\begin{math}
  w_{\underline{\kappa}}\colon r\mapsto w_{\underline{\kappa}}(r)
\end{math}
is a {\it subsolution\/}
for the backward initial value problem \eqref{back:FKPP:y(r)}.
It follows that
the backward initial value problem \eqref{back:FKPP:y(r)}
possesses a unique classical solution
$y\equiv y_{\mathrm{c}}(r)$ on the interval $(1 - \varrho, 1)$
which satisfies
\begin{math}
  y_{\mathrm{c}}(r)\geq w_{\underline{\kappa}}(r)
  = \underline{\kappa}\, (1-r)^{ 2 (\gamma_1 + \delta_1) }
\end{math}
for all $r\in (1 - \varrho, 1)$.
After returning to the original variables from eqs.\
\eqref{e:dU/dz} and \eqref{e:y=V^p'} we obtain
\begin{equation*}
  V(U)\geq \sqrt{\underline{\kappa}}\, (1-U)^{\gamma_1 + \delta_1}
  \quad\mbox{ for all }\, U\in (1 - \varrho, 1) \,.
\end{equation*}
We combine this inequality with the last limit in \eqref{e:d(1-)}
to conclude that there is a constant $c_2 > 0$ such that
\begin{equation}
\label{in:dU/dz}
{}- \genfrac{}{}{}0{\mathrm{d}z}{\mathrm{d}U}
  = \frac{d(U)}{V(U)} \leq
    \frac{c_2}{ (1-U)^{\gamma_1} }
  \quad\mbox{ for all }\, U\in (1 - \varrho, 1) \,.
\end{equation}
Notice that the relation
$(\gamma_1, \delta_1)\in \mathscr{M}_1^3$
implies also $\gamma_1 < 1$.
Consequently, fixing an arbitrary number
${\tilde U}\in (1 - \varrho, 1)$, denoting
${\tilde z} = z({\tilde U})\in (z_0,z_1)$,
and integrating ineq.~\eqref{in:dU/dz}
with respect to $U\in ({\tilde U}, 1)$, we arrive at
\begin{equation*}
    {\tilde z} - z_0
  = - \int_{\tilde U}^1
    \genfrac{}{}{}0{\mathrm{d}z}{\mathrm{d}U} \,\mathrm{d}U
  \leq c_2\int_{\tilde U}^1
    \frac{ \mathrm{d}U }{ (1-U)^{\gamma_1} }
  < \infty \,,
\end{equation*}
which forces $z_0 > -\infty$.

{\it Case 3:\/}$\quad$
$(\gamma_1, \delta_1)\in \mathscr{M}_1^2$.
This time we compare
$y\equiv y_{\mathrm{c}}\colon (0,1)\to (0,\infty)$
with the function
\begin{math}
  w_{\kappa}(r)\eqdef \kappa\, (1-r)^{\gamma_1 + \delta_1 + 1}
\end{math}
of $r\in [0,1]$, where $\kappa > 0$
is a suitable number to be determined later again.
From eq.~\eqref{oper:FKPP:w(r)} we deduce that there is
a sufficiently large number $\bar{\kappa} > 0$ such that
\begin{equation*}
  \mathscr{A}(w_{\bar{\kappa}})(r) \leq
  {}- \bar{\kappa} (\gamma_1 + \delta_1)\,
    (1-r)^{\gamma_1 + \delta_1} < 0
  \quad\mbox{ for all\/ $r\in (0,1)$ close to~$1$. }
\end{equation*}
Hence, there is a sufficiently small number $\varrho\in (0,1)$ such that
$w_{\bar{\kappa}}\colon r\mapsto w_{\bar{\kappa}}(r)$
is a {\it supersolution\/}
for the backward initial value problem \eqref{back:FKPP:y(r)}.

By similar arguments as above, we have
\begin{math}
  y_{\mathrm{c}}(r)\leq w_{\bar{\kappa}}(r)
  = \bar{\kappa}\, (1-r)^{\gamma_1 + \delta_1 + 1}
\end{math}
for all $r\in (1 - \varrho, 1)$.
After returning to the original variables from eqs.\
\eqref{e:dU/dz} and \eqref{e:y=V^p'} we obtain,
with a constant $c_3 > 0$,
\begin{equation}
\label{ine:dU/dz}
{}- \genfrac{}{}{}0{\mathrm{d}z}{\mathrm{d}U}
  = \frac{d(U)}{V(U)} \geq
    \frac{c_3}{ (1-U)^{ (\gamma_1 - \delta_1 + 1) / 2 } }
  \quad\mbox{ for all }\, U\in (1 - \varrho, 1) \,.
\end{equation}
Notice that the relation
$(\gamma_1, \delta_1)\in \mathscr{M}_1^2$
implies also
\begin{math}
  \frac{1}{2} (\gamma_1 - \delta_1 + 1)\geq 1 .
\end{math}
Again, we fix an arbitrary number
${\tilde U}\in (1 - \varrho, 1)$, denote
${\tilde z} = z({\tilde U})\in (z_0,z_1)$,
and integrate ineq.~\eqref{ine:dU/dz}
with respect to $U\in ({\tilde U}, 1)$, thus arriving at
\begin{equation*}
    {\tilde z} - z_0
  = - \int_{\tilde U}^1
    \genfrac{}{}{}0{\mathrm{d}z}{\mathrm{d}U} \,\mathrm{d}U
  \geq c_3\int_{\tilde U}^1
    \frac{ \mathrm{d}U }{ (1-U)^{ (\gamma_1 - \delta_1 + 1) / 2 } }
  = +\infty \,.
\end{equation*}
This estimate forces $z_0 = -\infty$.

{\it Case 4:\/}$\quad$
$(\gamma_1, \delta_1)\in \mathscr{M}_1^4$.
Finally, we compare
$y\equiv y_{\mathrm{c}}\colon (0,1)\to (0,\infty)$
with the function
\begin{math}
  w_{\kappa}(r)\eqdef \kappa\, (1-r)^{ 2 (\gamma_1 + \delta_1) }
\end{math}
of $r\in [0,1]$, where $\kappa > 0$
is a suitable number to be determined.
From eq.~\eqref{ope:FKPP:w(r)} we deduce that there is
a sufficiently large number $\bar{\kappa} > 0$ such that
\begin{equation*}
  \mathscr{A}(w_{\bar{\kappa}})(r) \leq
  {}- 2\bar{\kappa} (\gamma_1 + \delta_1)\,
  (1-r)^{ 2 (\gamma_1 + \delta_1) - 1 } < 0
  \quad\mbox{ for all\/ $r\in (0,1)$ close to~$1$. }
\end{equation*}
Hence, there is a sufficiently small number $\varrho\in (0,1)$ such that
$w_{\bar{\kappa}}\colon r\mapsto w_{\bar{\kappa}}(r)$
is a {\it supersolution\/}
for the backward initial value problem \eqref{back:FKPP:y(r)}.

Similarly as above, we have
\begin{math}
  y_{\mathrm{c}}(r)\leq w_{\bar{\kappa}}(r)
  = \bar{\kappa}\, (1-r)^{ 2 (\gamma_1 + \delta_1) }
\end{math}
for all $r\in (1 - \varrho, 1)$.
After returning to the original variables from eqs.\
\eqref{e:dU/dz} and \eqref{e:y=V^p'} we obtain,
with a constant $c_4 > 0$,
\begin{equation}
\label{neq:dU/dz}
{}- \genfrac{}{}{}0{\mathrm{d}z}{\mathrm{d}U}
  = \frac{d(U)}{V(U)} \geq
    \frac{c_4}{ (1-U)^{\gamma_1} }
  \quad\mbox{ for all }\, U\in (1 - \varrho, 1) \,.
\end{equation}
Notice that the relation
$(\gamma_1, \delta_1)\in \mathscr{M}_1^4$
implies also
\begin{math}
  \gamma_1\geq 1 .
\end{math}
Again, we fix an arbitrary number
${\tilde U}\in (1 - \varrho, 1)$, denote
${\tilde z} = z({\tilde U})\in (z_0,z_1)$,
and integrate ineq.~\eqref{neq:dU/dz}
with respect to $U\in ({\tilde U}, 1)$, thus arriving at
\begin{equation*}
    {\tilde z} - z_0
  = - \int_{\tilde U}^1
    \genfrac{}{}{}0{\mathrm{d}z}{\mathrm{d}U} \,\mathrm{d}U
  \geq c_4\int_{\tilde U}^1
    \frac{ \mathrm{d}U }{ (1-U)^{\gamma_1} }
  = +\infty \,,
\end{equation*}
which forces $z_0 = -\infty$.

The theorem is proved.
\qed
\par\vskip 10pt

\subsection{Comparisons with Previous Results}
\label{ss:Prev_Results}

The first result on the existence of travelling waves of
the so\--called {\em sharp\--type\/} for $c = c^{\ast}$
was obtained in
{\sc F.\ S\'anchez\--Gardu\~{n}o} and {\sc P.~K.\ Maini}
\cite[Theorem~2, p.~187]{Sanchez-Maini}.
The authors assume
$d(0) = 0$, $d > 0$ in $(0,1]$, $g(0) = g(1) = 0$, $g > 0$ in $(0,1)$,
and impose the following additional smoothness assumptions:
$d\in C^2([0,1])$, $d'(s) > 0$ and $d''(s)\neq 0$
for all $s\in [0,1]$,
$g\in C^2([0,1])$, $g'(0) > 0$ and $g'(1) < 0$.
These assumptions are weakened to
$d\in C([0,1])\cap C^1((0,1])$ and $g\in C([0,1])$ in
{\sc L.\ Malaguti} and {\sc C.\ Marcelli} \cite{Malag-Marce-1},
in Theorems 2, 3, and 14 (pp.\ 474, 475, and 493).
The authors in \cite{Malag-Marce-1} allow even for
$d'(0) = +\infty$ and $d(1) = 0$;
of particular interest to us are the existence results
for travelling waves of {\em sharp\--type\/}
\cite[Theorems 2(b) and 14(b)]{Malag-Marce-1}.

Our results are related to the existence results in
\cite{Malag-Marce-1}.
However, our results cover more general asymptotic behavior of
both terms, $d$ and $g$, near the equilibrium points $0$ and~$1$.
Indeed, their existence result
\cite[Theorem~2, p.~474]{Malag-Marce-1}
corresponds to the following parameter values in our case:
$\gamma_0 > 0$, $\delta_0 = 1$, $\gamma_1 > 0$, and $\delta_1 = 0$.
Another existence result in
\cite[Theorem~3, p.~475]{Malag-Marce-1} 
corresponds to our parameter values
$\gamma_0 + \delta_0 > 1$, $0 < \delta_0 < 1$, $\gamma_1 > 0$, and
$\delta_1 = 0$.
Furthermore, the existence result for doubly degenerate diffusion in
\cite[Theorem~14, p.~493]{Malag-Marce-1} corresponds to
$\gamma_0 > 0$, $\delta_0 = 1$, $\gamma_1 > 0$, and $\delta_1 = 1$.
In each of these cases, for $0 < \gamma_1 < 1$,
we obtain a wave profile $U$ with $z_0 > -\infty$, while for
$\gamma_1\geq 1$ we have $z_0 = -\infty$.

%
{\bf Acknowledgments.}
\begin{small}
The work of Pavel Dr\'abek was supported in part by
the Grant Agency of the Czech Republic (GA\v{C}R)
under Grant {\#}$18-03253$S.
\end{small}
%


\baselineskip=13.5pt
%
%
\makeatletter \renewcommand{\@biblabel}[1]{\hfill#1.} \makeatother
%
%


\begin{thebibliography}{99}

\bibitem{AronWein-1}
D.~G. Aronson and H.~F. Weinberger,
{\it Nonlinear diffusion in population genetics, combustion, and
     nerve propagation}, in:
{\sl ``Partial Differential Equations and Related Topics''},
     in {\em Lect. Notes in Math.}, Vol.~\textbf{446}, pp.\ 5--49.
     Springer\--Verlag, Berlin\--Heidelberg\--New York, 1975.

\bibitem{AronWein-2}
D.~G. Aronson and H.~F. Weinberger,
{\it Multi\-dimensional nonlinear diffusion arising in
     population genetics},
     Advances in Math., {\bf 30} (1978), 33--76.

\bibitem{DrabTak-1}
P. Dr\'abek and P. Tak\'a\v{c},
{\it New patterns of travelling waves in
     the generalized Fisher\--Kolmogorov equation},
     Nonlinear Differ. Equ. Appl. (NoDEA), {\bf 23}(2) (2016),
     Article~$7$ (online).
\color{blue}
     {\it Online}: http://dx.doi.org/10.1007/s00030-016-0365-2.
\color{black}

\bibitem{EGSanchez}
R. Engui\c{c}a, A. Gavioli, and L. Sanchez,
{\it A class of singular first order differential equations
     with applications in reaction\--diffusion},
     Discr. Cont. Dynam. Systems, {\bf 33}(1) (2013), 173--191.

\bibitem{Fife-McLeod}
P.~C. Fife and J.~B. Mc{L}eod,
{\it The approach of solutions of nonlinear diffusion equations to
     travelling front solutions},
     Arch. Rational Mech. Anal., {\bf 65}(4) (1977), 335--361.

\bibitem{Fisher}
R.~A. Fisher,
{\it The wave of advance of advantageous genes},
     Ann. of Eugenics, {\bf 7} (1937), 355--369.

\bibitem{Hamel-Nadira}
F. Hamel and N. Nadirashvili,
{\it Travelling fronts and entire solutions of
     the Fisher\--KPP equation in $\mathbb{R}^N$},
     Arch. Rational Mech. Anal., {\bf 157} (2001), 91--163.

\bibitem{KPP}
A.~N. Kolmogoroff, I.~G. Petrovsky, and N.~S. Piscounoff,
{\it \`Etude de l'\'equation de la diffusion avec croissance
     de la quantit\'e de mati\`ere et son application
     \`a un probl\`eme biologique},
     Bulletin Universit\'e d'\'Etat \`a Moscou
     (Bjul. Moskowskogo Gos. Univ.), S\'erie internationale,
     section~A, {\bf 1} (1937), 1--25.

\bibitem{Malag-Marce-1}
L. Malaguti and C. Marcelli,
{\it Sharp profiles in degenerate and doubly degenerate
     Fisher\--KPP equations},
     J.~Differential Equations, \textbf{195} (2003), 471--496.
\color{blue}
     {\it Online}: doi: 10.1016/j.jde.2003.06.005.
\color{black}

\bibitem{Murray}
J.~D. Murray,
{\sl ``Mathematical Biology''},
in {\em Biomathematics Texts}, Vol.~{\bf 19},
Springer-Verlag, Berlin--Hei\-del\-berg\--New York, 1993.

\bibitem{Murray-I}
J.~D. Murray,
{\sl ``Mathematical Biology I: An Introduction''}, $3$-rd Ed.
In {\em Interdisciplinary Applied Mathematics}, Vol.~{\bf 17},
Springer-Verlag, Berlin--Hei\-del\-berg\--New York, 2002.

\bibitem{Sanchez-Maini}
F. S\'anchez\--Gardu\~{n}o and P.~K. Maini,
{\it Existence and uniqueness of a sharp travelling wave
     in degenerate non\--linear diffusion Fisher\--PKK equation},
     J.~Math. Biology, \textbf{33} (1994), 163--192.

\bibitem{Tsoular_Wall}
A. Tsoularis and J. Wallace,
{\it Analysis of logistic growth models},
     Math. Biosciences, {\bf 179}(1) (2002), 21--55.

\end{thebibliography}
\end{document}